\documentclass[12pt,reqno]{amsart}

\usepackage{amssymb,latexsym}
\usepackage[matrix,arrow]{xy}
\usepackage[usenames]{color}
\usepackage{slashed}
\usepackage{parskip}
\newtheorem{thm}{Theorem}[section]
\newtheorem{cor}[thm]{Corollary}
\newtheorem{lem}[thm]{Lemma}

\theoremstyle{definition}

\theoremstyle{remark}
\newtheorem{rem}{Remark}[section]
\newcommand{\n}{\nabla}

\newcommand{\e}{\eqno}
\newcommand{\vi}{\varphi}
\newcommand{\s}{\Bbb{S}}

\newcommand{\ld}{\lambda}

\newcommand{\G}{\Gamma}
\newcommand{\om}{\omega}
\newcommand{\al}{\alpha}
\newcommand{\mr}{\mathfrak{R}}
\newcommand{\su}{\sum_{i=1}^k}
\newcommand{\de}{\Delta}
\newcommand{\sss}{\slashed{S}}

\newcommand{\ba}{\begin{aligned}}
\newcommand{\ea}{\end{aligned}}

\newcommand{\mg}{\mbox{grad}}
\begin{document}
\title{\text Extrinsic eigenvalue estimates of
the Dirac operator\bf}
\author []{ Daguang Chen }
% \address{\vspace{0.7cm}
\email { chendg@amss.ac.cn } \curraddr{Institute of Mathematics,
 Academy of Mathematics and Systems Sciences,
 Chinese Academy of Sciences, Beijing 100080, P.R. China.}
%\thanks { }
\date{29 Jan-2007}
%\begin{abstract} {}
%\end{abstract}
%\hspace{0.5cm}\begin{minipage}{12cm}
%\begin{center}
\begin{abstract}
For a compact spin manifold $M$ isometrically embedded into
Euclidean space, we derive the extrinsic estimates from above and
below for eigenvalues of the Dirac operators, which depend on the
second fundamental form of the embedding. We also show the bounds of
the ratio of the eigenvalues.

\end{abstract}
\maketitle
\renewcommand{\sectionmark}[1]{}
%-------------------------------------------------------------------------

\section{introduction}
The Laplacian operator and the Dirac operator are fundamental
differential operators. The estimates of their eigenvalues are
important in geometry, analysis and physics. Let $\Omega\subset\Bbb
R^n$ be the bounded domain of $\Bbb R^n$. Consider the eigenvalue
problem of the ~Laplacian
$$
\left \{\aligned & \Delta u=\xi u\qquad \mbox{in}\ \Omega,\\
&u|_{\partial \Omega}=0,\endaligned\right. \e(1.1)
$$
where $\Delta$ is the positive ~Laplacian in $\Bbb R^n$. It is
well known that this problem has a real and purely discrete
spectrum
$$
0<\xi_1<\xi_2\leq\xi_3\cdots\longrightarrow\infty.
$$
In 1956, using the technique of test-function,
 ~Payne, ~P\`{o}lya and
~Weinberger \cite{PPW}  proved,
$$
\xi_{k+1}-\xi_k\leq\frac{4}{nk}\su \xi_i.\e(1.2)
$$
In 1991, Yang \cite{Yang1,CY,CY1}  obtained very sharp inequality,
that is, he derived
$$
\su(\xi_{k+1}-\xi_i)(\xi_{k+1}-(1+\frac4n)\xi_i)\leq0.\e(1.3)
$$
From the inequality (1.3), we can obtain
$$
\xi_{k+1}\leq\frac1k(1+\frac4n)\su\xi_i.\e(1.4)
$$
The inequalities (1.3) and (1.4) are called Yang's first
inequality and second inequality, respectively. For the above
inequalities, we have the following relations
$$
(1.3)\Longrightarrow (1.4)\Longrightarrow(1.2).
$$
For the details, we refer to (\cite{Ash1,Ash2}) by M.S.Ashbaugh.
When $M$ is an n-dimensional compact manifold, there are the
similar results about the eigenvalue estimates for the ~Laplacian
(see, e.g.\cite{Li,CY}).

It is quite interesting to study the analogous estimates for the
Dirac operator. Let $M$ be a compact connected n-dimension
Riemannian manifold isometrically embedded into Euclidean space
$\Bbb R^N$ for some $N$. Estimates from above for the eigenvalues of
the Dirac operator can be attained in various ways(see
e.g.\cite{A,B,Baum,Bunke,Fr}). In 1991, N.Anhgel \cite{A} obtained
the analogous estimate of (1.2) for the Dirac operator. In the
present paper, by constructing suitable test ~spinors and
considering the corresponding Rayleigh quotients we obtain Yang-type
inequalities in term of the geometric data of ~Riemannian manifold.
Moreover, we obtain the lower bounds of eigenvalues for the square
of the Dirac operator. Finally, as an application of  \emph{theorem
4.1} in \cite{CY1}, a combination of \emph{Theorem 3.4} with
\emph{Theorem 4.2} yields the bounds of the ratio $\mu_{k+1}/\mu_1$,
where $\mu_i=\ld_i+\frac14\underset{M}\sup h^2$, $\{\ld_i\}_{1\leq
i\leq\infty}$ are the eigenvalues for the square of the Dirac
operator.
%---------------------------------------------------------------------

\section{Preliminaries}

Let $M$ be an n-dimensional compact ~Riemannian manifold and let
$x=(x^1,\cdots,x^N):M\hookrightarrow\mathbb{ R}^N$ be an isometric
embedding of $M$ in $\mathbb {R}^N$. For a fixed point $x\in M$ and
an orthonomal basis $\{e_i,e_\alpha\}$ with $e_i$ tangent to $M$ and
$e_\alpha$ normal to $M$ for $1\leq i\leq n, n+1\leq \alpha \leq n$.
Then we have the structure equation
$$
\begin{aligned}
&dx=\om^ie_i,\qquad \om^\alpha=0,\\
&de_i=\om^j_ie_j+\om^\alpha_ie_\alpha,
\qquad \om^\alpha_i=h^\alpha_{ij}\om^j,\\
&de_\alpha=\om^j_\alpha e_j+\om^\beta_\alpha e_\beta,
\end{aligned}\eqno(2.1)
$$
where $h^\alpha_{ij}$ is components of the second fundamental
form, $\omega^i_j$ is the connection form,
 $\omega^\beta_\alpha$ is the normal connection form.  \\
 One can infer that, ~pointwise on $M$
 $$
\sum_{ p=1}^N|\mg x^p|^2=n.\e(2.2)
 $$
Let $dx=x_i\om^i$, comparing the first line of (2.1), we have
$x_i=e_i$. From the structure equation (2.1), we have
$$\begin{aligned}
x_{ij} \om^j&=dx_i-x_j\om_i^j\\
&=de_i-\om_i^je_j\\
&=\om_i^je_j+\om_i^\alpha e_\alpha-\om_i^je_j\\
&=h^\alpha_{ij}e_\alpha\om^j.
\end{aligned}
$$
Therefore
$$
\Delta x=-\sum h_{ii}^\al e_\al=-\mathbf{H}.\e(2.3)
$$
where $\mathbf{H}=\sum h^\al_{ii}e_\al$ is the mean curvature
vector field of $M$, $\Delta$ is the positive ~Laplacian. We
denote by $H=\frac 1n|\mathbf H|$ the mean curvature function.
Since the ambient space is Euclidean space, we have
$$
S=n^2H^2-h^2\e(2.4)
$$
where $h^2=\sum_{i,j,\alpha}(h^\alpha_{ij})^2$, $S$ is the scalar
curvature of $M$.

Now we review the following basic facts in \cite{LM} or in
\cite{Fr}. Let $\sss$ be a Dirac bundle over the ~Riemannian
manifold $M$. This means that $\sss$ is a bundle of left modules
over $Cl(M)$ together with a ~Riemannian metric and connection on
$\sss$ having the following properties:
 \begin{itemize}
\item $<e\sigma_1,\sigma_2>+<\sigma_1,e\sigma_2>=0, $ \item
$\n(\vi\sigma)=\n\vi\sigma+\vi\n\sigma,$
\end{itemize}
for
$\sigma_1,\sigma_2,\sigma\in\G(\sss),\vi\in\G(Cl(M)),e\in\G(TM)$,
where $Cl(M)$ is the Clifford bundle over $M$, the covariant
derivative on $\sss$ also denote by $\n$.

The Dirac operator is a first order elliptic differential operator
$D:\G(\sss)\longrightarrow \G(\sss)$, which is locally given by
$$D=\textstyle\sum_{i=1}^n e_i\n_i.$$
For $f\in C^{\infty}(M)$ and $\sigma\in\G(\sss)$, we get
$$
D(f\sigma)=\mg( f) \sigma+fD \sigma\e(2.5)
$$
$$
D^2(f \sigma)=\Delta (f)\sigma-2\n_{\mg
(f)}\sigma+fD^2\sigma,\e(2.6)
$$
where $\Delta $ is the positive scalar ~Laplacian. For the Dirac
operator $D$, one has the general ~Bochner identity
$$D^2=\n^*\n+\mathfrak{R},\e(2.7)$$
where $\mathfrak{R}$ is the curvature morphism acting on Dirac
bundle $\sss$, which is given by
$$\ba&\mathfrak{R}=\sum_{i<j}e_ie_j\mathcal
{R}_{e_i,e_j},\\
&\mathcal{R}_{e_i,e_j}=[\n_i,\n_j]-\n_{[e_i,e_j]}.\ea$$

There exist two basic cases of the Dirac bundles:
\begin{itemize}
\item The Clifford bundle $Cl(M)$. The Dirac operator in this case
is a square root of the classical ~Hodge ~Laplacian.\item The
~spinor bundles. Suppose $M$ is spin manifold with the spin
structure on its tangent bundle. Let $\s$ be any ~spinor bundle
associated to tangent bundle. Then $\s$ carries a canonical
~Riemannian connection. In this case, $D$ is the classical Dirac
operator (also ~Atiyah-Singer operator) and $\mathfrak{R}=\frac 14
S$, where $S$ is the scalar curvature of $M$.
\end{itemize}

%------------------------------------------------------------------
\section{upper bounds}
It is well known that the Dirac operator $D$ of any Dirac bundle
is self-adjoint and elliptic on the compact manifold. Therefore,
$D^2$ has a discrete spectrum contained in $\Bbb R$ numbered like
$$
0\leq\ld_1\leq\ld_2\leq\cdots\nearrow\infty.
$$
and one can find an orthonormal basis $\{\vi_j\}_{j\in\Bbb N}$ of
$L^2(\sss)$ consisting of  eigenfunctions of $D^2$
(i.e.\,$D^2\vi_j=\ld_j\vi_j, j\in \Bbb N$).  Such a system
$\{\ld_j;\vi_j\}_{j\in\Bbb N}$ is called a $spectral$
$decomposition$ of $L^2(\sss)$ generated by $D^2$, or, in short, a
$spectral$ $resolution$ of $D^2$. In the following, we  want to
estimate the eigenvalue $\ld_{k+1}$ in terms of the previous ones,
$\ld_1,\ld_2,\cdots,\ld_k$, their eigenspinors
$\vi_1,\cdots,\vi_k$, the immersion $x$, and the curvature term
$\mathfrak{R}$.

In the present paper, we denote $(\cdot ,\cdot)=\Re\int_M< \cdot
,\cdot >$ on Dirac bundle $\sss$.
%---------------------------------------------------------------------------------------
\begin{thm}
Let $M$  be a compact ~Riemannian manifold of dimension $n$ and
$x=(x^1,\cdots,x^N):M\longrightarrow \Bbb{R}^N$ be an isometric
embedding with mean curvature function $H$. Let $D$ be the Dirac
operator of any Dirac bundle $\sss$ over $M$ and let
$\{\ld_j;\vi_j\}_{j\in\Bbb N}$ be a spectral resolution of $D^2$.
Then
$$
\sum_{i=1}^{k}(\ld_{k+1}-\ld_i)^2\leq\frac4n
\sum_{i=1}^k(\ld_{k+1}-\ld_i)\left(\ld_i+\frac14\Big(
n^2(H^2\vi_i,\vi_i)-4\big(\mr \vi_i,\vi_i\big)\Big)\right)\e(3.1)
$$
Moreover, if $(\mr\vi_i,\vi_i)\geq R$ for some $R\in\Bbb R$,
denoting $\eta_i=\ld_i+\frac14 n^2\underset{M}\sup{H^2}-R$, then
we have the following Yang-type inequality
$$
\sum_{i=1}^{k}(\eta_{k+1}-\eta_i)\Big(\eta_{k+1}-(1+\frac4n)\eta_i\Big)\leq
0.\e(3.2)
$$
\end{thm}
%---------------------------------------------------------------------------------------

\begin{proof}
Since $\{\ld_j;\vi_j\}_{j\in\Bbb N}$ is a spectral resolution of
$D^2$, $\int<\vi_i,\vi_j>=\delta_{ij} $ implies
$(\vi_i,\vi_j)=\delta_{ij}$, for $ \ \forall\  i,j\in\Bbb N. $
 Letting $g=x^p,1\leq p\leq N$ and for
$1\leq i,j\leq k$, we take
$$
\left \{ \aligned
a_{ij}&=(g\vi_i,\vi_j),\\
\psi_i&=g\vi_i-\su a_{ij}\vi_j,\\
b_{ij}&=(\vi_i,\frac12\Delta g\vi_j-\nabla_{\mg(g)}\vi_j).
\endaligned \right .\e(3.3)
$$
Then one can infer $a_{ij}=a_{ji}$,  $ (\psi_i,\vi_j)=0, \quad \text{for} \quad j=1,\cdots, k$.\\
According to Green's formula and (2.6), we infer
$$
\begin{aligned}
\ld_i a_{ij}&=(gD^2\vi_i,\vi_j)\\
&=(\vi_i,D^2(g \vi_j))\\
&=(\vi_i,\Delta g\vi_j-2\nabla_{\mg(g)}\vi_j+gD^2\vi_j)\\
&=\ld_ja_{ij}+2(\vi_i,\frac12\Delta g \vi_j-\n_{\mg(g)}\vi_j)\\
&=\ld_ja_{ij}+2b_{ij}.
\end{aligned}
$$
Hence, we get
$$
2b_{ij}=(\ld_i-\ld_j)a_{ij}=-2b_{ji}.\e(3.4)
$$
From (2.6) and (3.3), one immediately has
$$ D^2\psi_i=\de
g\vi_i-2\n_{\mg(g)}\vi_i+\ld_ig\vi_i-\sum_{j=1}^k\ld_ja_{ij}\vi_j,
$$
Therefore
$$
\ba (D^2\psi_i,\psi_i)&=(\de
g\vi_i-2\n_{\mg(g)}\vi_i+\ld_ig\vi_i-\sum_{j=1}^k\ld_ja_{ij}\vi_j,\psi_i)\\
&=(\de g\vi_i-2\n_{\mg(g)}\vi_i,\psi_i)+\ld_i\|\psi_i\|^2.
\ea\e(3.5)
$$
From the Rayleigh quotient, one gets
$$
\ld_{k+1}\leq \frac{(D^2\psi_i,\psi_i)}{\|\psi_i\|^2}.\e(3.6)
$$
By (3.5), (3.6) can be written as
$$
(\ld_{k+1}-\ld_i)\|\psi_i\|^2\leq (\de
g\vi_i-2\n_{\mg(g)}\vi_i,\psi_i).\e(3.7)
$$
But (2.6) yields
$$
\ba
 (\de g\vi_i-2\n_{\mg(g)}\vi_i,\psi_i)&=(\de
 g\vi_i-2\n_{\mg(g)}\vi_i,g\vi_i)-
 \sum_{j=1}^ka_{ij}\left(D^2(g\vi_i)-gD^2\vi_i,\vi_j\right)\\
 &=(\de
 g\vi_i-2\n_{\mg(g)}\vi_i,g\vi_i)+\sum_{j=1}^k(\ld_i-\ld_j)a_{ij}^2.
\ea
$$
Therefore, we obtain
$$
(\ld_{k+1}-\ld_i)\|\psi_i\|^2\leq(\de
g\vi_i-2\n_{\mg(g)}\vi_i,g\vi_i)+\sum_{j=1}^k(\ld_i-\ld_j)a_{ij}^2.\e(3.8)
$$
From ~Schwarz inequality and (3.7), we obtain
$$
\ba(\ld_{k+1}&-\ld_i)\left[\left(\de
g\vi_i-2\n_{\mg(g)}\vi_i,\psi_i\right)\right]^2\\
=&(\ld_{k+1}-\ld_i)\left[\left(\de
g\vi_i-2\n_{\mg(g)}\vi_i+2\sum_{j=1}^k
b_{ij}\vi_j,\psi_i\right)\right]^2\\
\leq & (\ld_{k+1}-\ld_i)\|\psi_i\|^2\|\de
g\vi_i-2\n_{\mg(g)}\vi_i+2\sum_{j=1}^k
b_{ij}\vi_j\|^2\\
\leq &\left(\de g\vi_i-2\n_{\mg(g)}\vi_i,\psi_i\right) \|\de
g\vi_i-2\n_{\mg(g)}\vi_i+2\sum_{j=1}^k b_{ij}\vi_j\|^2. \ea
\e(3.9)
$$
Multiplying (3.9) by $(\ld_{k+1}-\ld_i)$ and taking sum on $i$
from 1 to k, we have
$$\ba \su (\ld_{k+1}-\ld_i)^2&\left(\de
g\vi_i-2\n_{\mg(g)}\vi_i,g\vi_i\right)+\sum_{i,j=1}^k
(\ld_i-\ld_j)(\ld_{k+1}-\ld_i)^2a_{ij}^2\\
&\leq \su (\ld_{k+1}-\ld_i) \|\de
g\vi_i-2\n_{\mg(g)}\vi_i+2\sum_{j=1}^k b_{ij}\vi_j\|^2.
\ea\e(3.10)
$$
From the Eq. (2.4) in \cite{A}, we have
$$
\left(\de g\vi_i-2\n_{\mg (g)}\vi_i,g\vi_i\right)=\int |\mg
g|^2|\vi_i|^2.\e(3.11)
$$
Therefore (2.2) and (3.11) yield
$$
 \sum_{p=1}^{N}\su(\ld_{k+1}-\ld_i)^2
 \left(\de g\vi_i-2\n_{\mg (g)}\vi_i,g\vi_i\right)
=n\su (\ld_{k+1}-\ld_i)^2.\e(3.12)
$$

In order to complete the proof of theorem, we need the following
Lemma:
\begin{lem} With the above notations, then
$$\ba
&\sum_{p=1}^{N}\|\de g\vi_i-2\n_{\mg(g)}\vi_i+2\sum_{j=1}^k b_{ij}\vi_j\|^2\\
&=n^2(H^2\vi_i,\vi_i)+4\ld_i-4(\mr
\vi_i,\vi_i)-4\sum_{p=1}^{N}\sum_{j=1}^kb_{ij}^2. \ea\e(3.13)
$$
\end{lem}
\begin{proof} By the \emph{Lemma 2.9} in \cite{A}, a straightforward
computation yields
$$\ba
&\sum_{p=1}^{N}\|\de g\vi_i-2\n_{\mg(g)}\vi_i
+2\sum_{j=1}^k b_{ij}\vi_j\|^2\\
& =\sum_{p=1}^{N}\|\de
g\vi_i-2\n_{\mg(g)}\vi_i\|^2+4\sum_{p=1}^{N}\sum_{j=1}^kb_{ij}^2
+4\sum_{p=1}^{N}\sum_{j=1}^kb_{ij}\left(\de
g\vi_i-2\n_{\mg(g)}\vi_i,\vi_j\right)\\
&=n^2(H^2\vi_i,\vi_ i)+4\ld_i-4(\mr
\vi_i,\vi_i)+8\sum_{p=1}^{N}\sum_{j=1}^kb_{ij}b_{ji}+4\sum_{p=1}^{N}\sum_{j=1}^kb_{ij}^2\\
&=n^2(H^2\vi_i,\vi_i)+4\ld_i-4(\mr
\vi_i,\vi_i)-4\sum_{p=1}^{N}\sum_{j=1}^kb_{ij}^2. \ea
$$
\end{proof}

From (3.10), (3.12) and (3.13), one gets
$$\ba
n\su
(\ld_{k+1}-\ld_i)^2+&\sum_{p=1}^{N}\sum_{i,j}^k(\ld_i-\ld_j)(\ld_{k+1}-\ld_i)^2a_{ij}^2\\
&\leq\su
(\ld_{k+1}-\ld_i)\left[n^2(H^2\vi_i,\vi_i)-4(\mr\vi_i,\vi_i)\right]\\&+4\su
(\ld_{k+1}-\ld_i)\ld_i
-4\sum_{p=1}^{N}\sum_{i,j=1}^k(\ld_{k+1}-\ld_i)b_{ij}^2.
\ea\e(3.14)
$$
From the relation (3.4), one can infer
$$
\sum_{i,j}^k(\ld_i-\ld_j)(\ld_{k+1}-\ld_i)^2a_{ij}^2=
-\sum_{i,j}^k(\ld_{k+1}-\ld_i)(\ld_i-\ld_j)^2a_{ij}^2=
-4\sum_{i,j=1}^k(\ld_{k+1}-\ld_i)b_{ij}^2.\e(3.15)
$$

Finally, we obtain
$$\ba
n\su (\ld_{k+1}-\ld_i)^2 \leq&\su
(\ld_{k+1}-\ld_i)\left[n^2(H^2\vi_i,\vi_i)-4(\mr\vi_i,\vi_i)\right]\\&+4\su
(\ld_{k+1}-\ld_i)\ld_i. \ea
$$
\end{proof}
\begin{rem}
\emph{Theorem 3.1} improve \emph{Theorem 3.1} in \cite{A}, that
is,
$$
\ld_{k+1}-\ld_k\leq \frac
nk\su(H^2\vi_i,\vi_i)+\frac{4}{nk}\su\ld_i-\frac
4{nk}\su(\mr\vi_i,\vi_i).
$$

\end{rem}
From the \emph{Theorem 3.1}, we also obtain Yang's second
inequality.
\begin{cor} Under the same assumptions in \emph{Theorem 3.1}, then
one has
$$
\eta_{k+1}\leq(1+\frac4n)\frac1k\su \eta_i.\e(3.16)
$$
\end{cor}
From (2.4) and \emph{Theorem 3.1}, we immediately deduce
\begin{thm} Let $M$ be an n-dimensional compact ~Riemannian spin
manifold and $x:M\longrightarrow\Bbb R^N$ be isometrically embedded
into $\Bbb R^N$. Let $(\ld_i,\vi_i)_{i\in\Bbb N}$ be the spectral
resolution of the square of the Dirac operator $D$, i.e.
$D^2\vi_i=\ld_i\vi_i$. Order the eigenvalues in an increasing
sequence
$$
0\leq\ld_1\leq\ld_2\leq\cdots\nearrow\infty.
$$
Then one gets the Yang's first inequality
 $$
\su\Big(\mu_{k+1}-\mu_i\Big)\left(\mu_{k+1}-\Big(1+\frac4n\Big)\mu_i\right)\leq
0 \e(3.17)
$$
and the Yang's second inequality
$$
\mu_{k+1}\leq\left(1+\frac4n\right)\frac1k\su \mu_i,\e(3.18)
$$
where $\mu_i=\ld_i+\frac14(n^2\sup H^2-\inf S)$, $H, S$ are the
mean curvature function and scalar curvature of $M$ respectively.
\end{thm}

%-------------------------------------------------------------------
 \begin{rem} From  (3.18), one can deduce  the
 \emph{Theorem 3.7} in \cite{A}.
That is, one can infer
$$
\mu_{k+1}-\mu_k\leq \frac4{nk}\su \mu_i.
$$
\end{rem}
%-------------------------------------------------------------------

If the isometric embedding $M\xrightarrow[]{x} \Bbb
S^{N-1}(1)\xrightarrow[]{i}\Bbb R^N$ is the minimal embedding into
the Euclidean unit sphere $\Bbb S^{N-1}\subset\Bbb R^N$, then the
coordinate functions $x^1,x^2,\cdots,x^N$ are the eigenfunctions of
the ~Laplacian with the eigenvalue $n$, and $\mathbf
H=-n(x^1,x^2,\cdots,x^N)$ and the mean curvature function $H\equiv
1$. Thus in the case of minimally embeded submanifolds into
Euclidean spheres, one gets Yang-type inequalities
\begin{cor}
Let $M$  be a compact ~Riemannian spin manifold of dimension $n$ and
$x=(x^1,\cdots,x^N):M\xrightarrow[]{x} \Bbb
S^{N-1}(1)\xrightarrow[]{i}\Bbb R^N$ be the minimal embedding into
the Euclidean unit sphere $\Bbb S^{N-1}\subset\Bbb R^N$. Let $D$ be
the Dirac operator of spinor bundle $\sss$ over $M$ and let
$\{\ld_j;\vi_j\}_{j\in\Bbb N}$ be a spectral resolution of $D^2$.
Then one gets Yang's first inequality
$$
\su\Big(\nu_{k+1}-\nu_i\Big)\left(\nu_{k+1}-
\Big(1+\frac4n\Big)\nu_i\right)\leq 0 \e(3.19)
$$
and Yang's second inequality
$$
 \nu_{k+1}\leq\left(1+\frac4n\right)\frac1k\su \nu_i,\e(3.20)
$$
where $\nu_i=\ld_i+\frac14 (n^2-\underset{M}\inf S)$, $S$ is
scalar curvature of $M$.
\end{cor}

%------------------------------------------------------------------------
\section{lower bounds}
 In this section, we'll show lower bounds of
 the eigenvalues for the Dirac operator.
%------------------------------------------------------------------------
\begin{thm} Under the same assumptions and the notations
 as \emph{Theorem 3.1}, one gets
$$
\sum_{i=1}^n(\eta_{i+1}-\eta_1)\leq 4\eta_1.\e(4.1)
$$
\end{thm}
%-------------------------------------------------------------------------
\begin{proof}Since $\{\ld_j;\vi_j\}_{j\in\Bbb N}$ is a spectral resolution of
$D^2$, $\int<\vi_i,\vi_j>=\delta_{ij} $ implies
$(\vi_i,\vi_j)=\delta_{ij}$, for $\forall\  i,j\in\Bbb Z$.  First
define the matrix $B=\Big((x^p\vi_1,\vi_{q+1})\Big)_{1\leq p,q\leq
N}$. From the \textbf{QR}-factorization theorem, one gets $R=QB$,
where $Q=(Q_q^p)$ is an orthogonal $N\times N$ matrix and $R$ is
real upper triangular matrix. That is,
$$
(g^p\vi_1,\vi_{r+1})=(Q_q^px^q\vi_1,\vi_{r+1})=Q_q^p(x^q\vi_1,\vi_{r+1})
 =0,\qquad \mbox{for} \quad 1\leq r<p\leq N.\e(4.2)
$$
where $g^p=Q_q^px^q$. Second define the test spinors
$$
\Psi_p=g^p\vi_1-(g^p\vi_1,\vi_1)\vi_1.\e(4.3)
$$
Then
$$
(\Psi_p,\vi_{r+1})=0,\qquad \mbox{for}\ 0\leq r<p.
$$
On one hand, from Rayleigh quotient, we have
$$
\ld_{p+1}\leq\frac{\|D^2\Psi_p\|^2}{\|\Psi_p\|^2}.\e(4.4)
$$
 One the other hand, from (2.6), one gets
 $$
 D^2\Psi_p=\Delta g^p\vi_i-2\n_{\mg
 (g^p) }\vi_1+\ld_1(g^p\vi_1,\vi_1).\e(4.5)
$$
Therefore, (4.4) can be written as
$$
(\ld_{p+1}-\ld_1)\|\Psi_p\|^2\leq (\Delta
g^p\vi_1-2\n_{\mg{(g^p)}}\vi_1,\Psi_p).\e(4.6)
$$
From ~Schwarz inequality, we obtain
$$
\left(\Delta g^p\vi_1-2\n_{\mg{(g^p)}}\vi_1,\Psi_p\right)^2\leq
\|\Psi_p\|^2\|\Delta g^p\vi_1-2\n_{\mg{(g^p)}}\vi_1\|^2.\e(4.7)
$$
Multiplying (4.7) by $(\ld_{p+1}-\ld_1)$ and then by (4.6), one
gets
$$
\ba(\ld_{p+1}-\ld_1)&(\Delta
g^p\vi_1-2\n_{\mg{(g^p)}}\vi_1,\Psi_p)^2\\ &\leq
(\ld_{p+1}-\ld_1)\|\Psi_p\|^2\|\Delta
g^p\vi_1-2\n_{\mg{(g^p)}}\vi_1\|^2\\
&\leq (\Delta g^p\vi_1-2\n_{\mg{(g^p)}}\vi_1,\Psi_p)\|\Delta
g^p\vi_1-2\n_{\mg{(g^p)}}\vi_1\|^2. \ea
$$
That is,
$$
(\ld_{p+1}-\ld_1)(\Delta g^p\vi_1-2\n_{\mg{g^p}}\vi_1,\Psi_p)\leq
 \|\Delta g^p\vi_1-2\n_{\mg(g^p)}\vi_1\|^2.\e(4.8)
$$
By \emph{Lemma 2.9} in \cite{A} and by the definition of $g^p$, we
have
$$
 \sum_{p=1}^N\|\Delta
 g^p\vi_1-2\n_{\mg{(g^p)}}\vi_1\|^2= 4\ld_1
 +n^2(H^2\vi_1,\vi_1)-4(\mr \vi_1,\vi_1).\e(4.9)
$$
From (2.4) in \cite{A}, one obtains
$$
(\Delta
g^p\vi_1-2\n_{\mg(g^p)}\vi_1,\Psi_p)=
 (\Delta
g^p\vi_1-2\n_{\mg(g^p)}\vi_1,g^p\vi_1)=
(|\mg(g^p)|^2\vi_1,\vi_1).\e(4.10)
$$
Then we deduce
$$
\sum_{p=1}^N(\ld_{p+1}-\ld_1)(|\mg(g^p)|^2\vi_1,\vi_1)\leq 4\ld_1
 +n^2(H^2\vi_1,\vi_1)-4(\mr \vi_1,\vi_1).\e(4.11)
$$
From (2.2) and $Q\in O(N)$, one infers
$$
\ba \sum_{p=1}^N|\mg(g^p)|^2&=\sum_{p=1}^N\sum_{i=1}^n|\n_ig^p|^2
=\sum_{p=1}^N\sum_{i=1}^n\sum_{q,r=1}^NQ^p_qQ^p_r(\n_ix^q,\n_ix^r)\\
&=\sum_{q,r=1}^N\sum_{i=1}^n\delta_{qr}(\n_ix^q,\n_ix^r)\\
&=\sum_{p=1}^N|\mg(x^p)|^2=n. \ea$$

Therefore, (4.11) becomes
$$
\sum_{p=1}^N\ld_{p+1}(|\mg(g^p)|^2\vi_1,\vi_1)\leq (4+n)\ld_1
 +n^2(H^2\vi_1,\vi_1)-4(\mr \vi_1,\vi_1).\e(4.12)
$$
Let $P$ be a point of $M$. It is possible to take a coordinate
system $(\bar x^1,\cdots,\bar x^N)$  with origin $ P$, $x=x(P)+\bar
x A$, $A\in O(N)$ such that
 $T_PM=\text{span}\{(\frac {\partial }{\partial \bar x^1})_P,\cdots,
 (\frac {\partial }{\partial \bar x^n})_P\}$. In fact, we have
$$\n\bar x^{n+1}=\cdots=\n\bar x^{N}=0,\qquad\n_i\bar
x^j=\delta_i^j, \qquad(i,j=1,\cdots,n).\e(4.13)
  $$
From (4.13) and the  matrices $Q,A\in O(N)$, one gets
$$\ba
  |\mg(g^p)|^2&=\sum_{i=1}^n\left|\textstyle
  \sum_{q,r=1}^NQ^p_qa^q_r\n_i\bar x^r\right|^2\\
&=\sum_{i=1}^n\left|\textstyle
  \sum_{q=1}^NQ^p_qa^q_i\right|^2\\
  &\leq\sum_{r=1}^N\left|\textstyle
  \sum_{q=1}^NQ^p_qa^q_r\right|^2\\
  &=1,
\ea$$
since $\left(\textstyle
  \sum_{q=1}^NQ^p_qa^q_r\right)_{1\leq p,r\leq N}$
  is also an orthogonal
  matrix.

Thus we have
 $$\ba
  \sum_{p=1}^N\ld_{p+1}|\mg(g^p)|^2&
  \geq \sum_{i=1}^n\ld_{i+1}|\mg(g^i)|^2+
 \ld_{n+1} \sum_{\alpha=n+1}^N|\mg(g^\alpha)|^2\\&
 =\sum_{i=1}^n\ld_{i+1}|\mg(g^i)|^2+
 \ld_{n+1} (n-\sum_{i=1}^n|\mg(g^i)|^2)\\&
 =\sum_{i=1}^n\ld_{i+1}|\mg(g^i)|^2+
 \ld_{n+1} \sum_{i=1}^n(1-|\mg(g^i)|^2)\\&
 \geq \sum_{i=1}^n\ld_{i+1}|\mg(g^i)|^2+
  \sum_{i=1}^n\ld_{i+1}(1-|\mg(g^i)|^2)\\&
  =\sum_{i=1}^n\ld_{i+1}.
\ea\e(4.14) $$
 A combination (4.12) with (4.14) yields Theorem 4.1.
\end{proof}
Obviously,\emph{Theorem 4.1} implies
\begin{thm}
Under the same assumptions and the notations as
 \emph{Theorem 3.4}, then one gets
$$
\sum_{i=1}^n(\mu_{i+1}-\mu_1)\leq 4\mu_1. \e(4.15)
$$
\end{thm}

\begin{cor}
Under the same assumptions and the notations as \emph{Corollary
3.5}, then
$$
\su(\nu_{k+1}-\nu_1)\leq 4\nu_1.\e(4.16)
$$
\end{cor}
\begin{rem}
For the eigenvalue problem (1.1), there exists the
inequality\cite{PPW}
$$
 \su(\xi_{i+1}-\xi_i)\leq 4\xi_1
$$
 where $\xi_i$ is the i-st eigenvalue of the ~Laplacian.
\end{rem}

 From the \emph{Theorem 3.1} in \cite{CY1} , a combination
 \emph{Theorem 3.4}  with  \emph{Theorem 4.1}   yields:
\begin{thm}Under the same assumptions  and notations as \emph{Theorem 3.4},
 then\\
 (1) for $n\geq 41$ and $k\geq 41,$
  $$
  \mu_{k+1}\leq k^{2/n}\mu_1;
  $$
 (2) for any n and k,
 $$
\mu_{k+1}\leq \left(1+\frac {a(min\{n,k-1\})}{n}\right)k^{2/n}\mu_1,
 $$
where the bound of $a(m)$ can be formulated as:
$$
\left \{ \aligned
 a(1)&\leq 2.64,\\
 a(m)&\leq 2.2-4\log\left(1+\frac1{50}(m-3)\right),\qquad \mbox{for}\quad
 m\geq 2.
\endaligned \right .
$$
\end{thm}

%----------------------------------------------------------------------------
From the \emph{Corollary 2.1} in \cite{CY1} or \cite{CY}, we can
deduce the simple and clear inequality
$$
\mu_{k+1}\leq\left(1+\frac4n\right)k^{\frac2n}\mu_1  \e(4.17)
$$
where $\mu_i$ is given in  \emph{Theorem 3.4}. From the
\emph{Lemma 1.12.6} in \cite{Gilkey}, (4.17) is a best possible
estimate of $\mu_{k+1}$ in the sense of order.

%-----------------------------------------------------------------------------------------------
\textbf{Acknowledgments} The author would like to  express his
gratitude to Professor X. \-Zhang and H.C. \-Yang for their
encouragements, suggestions and support.

%-----------------------------------------------------------------

\providecommand{\bysame}{\leavevmode\hbox
to3em{\hrulefill}\thinspace}

\end{document}